\documentclass[12pt]{amsart}
\usepackage{geometry}                
\usepackage{graphicx}
\usepackage{amssymb}
\usepackage{epstopdf}
\usepackage{xcolor,color}
\usepackage{amsfonts}
\usepackage{amsmath,amssymb}
\newtheorem{proposition}{Proposition}
\newtheorem{definition}{Definition}
\newtheorem{theorem}{Theorem}

\title[A note on the periodic orbits of Wolbachia spread dynamics]{A note on the periodic orbits of Wolbachia spread dynamics in mosquito populations in periodic environments}
\author{J. S. C\'{a}novas }

\begin{document}

\begin{abstract}
We consider the periodic model introduced in \cite{zheng}, and disprove the conjectures on the number of periodic orbits the model can have. We rebuild the conjecture to prove that for periodic sequences of maps of any period, the number of non-zero periodic trajectories is bounded by two.


\end{abstract}
\maketitle

\newtheorem{lemma}[theorem]{Lemma}
\newtheorem{conjecture}{Conjecture}
\newtheorem{example}[theorem]{Example}
\newtheorem{corollary}[theorem]{Corollary}
\newtheorem{remark}[theorem]{Remark}
\newtheorem{examples}[theorem]{Examples}
\newtheorem{remarks}[theorem]{Remarks}
\newtheorem{notation}[theorem]{Notation}


\section{Introduction}

In the paper \cite{zheng}, the authors introduced and studied a model of Wolbachia
spread dynamics in mosquito populations in cyclic environments. The model is
derived from the autonomous model given by
\begin{equation*}
x_{n+1}=\frac{(1-\mu )(1-s_{f})x_{n}}{s_{h}x_{n}^{2}-(s_{h}+s_{f})x_{n}+1},
\end{equation*}%
where $x_{n}\in \lbrack 0,1]$ is the frequency of Wolbachia-infected adults
at the $n$-th generation, and $\mu \in \lbrack 0,1)$, $s_{f}\in \lbrack 0,1)$
and $s_{h}\in (0,1]$. This model was introduced in \cite{fine} by modifying
a simpler model introduced in \cite{caspari}, which can be obtained taking the parameter $\mu=0$. We refer the reader to \cite{caspari,fine,zheng} and the references therein for more biological knowledge of this
model.

Let%
\begin{equation}
f(x)=\frac{(1-\mu )(1-s_{f})x}{s_{h}x^{2}-(s_{h}+s_{f})x+1}.  \label{map}
\end{equation}%
The following facts are easily checked.

\begin{enumerate}
\item[(F1)] The denominator of the map $f(x)$ vanishes for%
\begin{equation*}
o_{\pm }=\frac{s_{h}+s_{f}\pm \sqrt{(s_{h}+s_{f})^{2}-4s_{h}}}{2s_{h}}
\end{equation*}%
which are real numbers if and only if%
\begin{equation*}
s_{f}\geq 2\sqrt{s_{h}}-s_{h}.
\end{equation*}%
In this case, as $s_{f}\leq 1$, we see that $o_{-}\geq 1$ and hence the map $%
f$ is properly defined in $[0,1]$. Besides, it is straightforward to check that $(s_{h}+s_{f})^{2}-4s_{h}<0$ whenever $s_f<s_h$, and so no real poles exist.

\item[(F2)] The derivative of $f$ is%
\begin{equation*}
f^{\prime }(x)=-\frac{(\mu -1)(s_{f}-1)(s_{h}x^{2}-1)}{%
(s_{h}x^{2}-(s_{h}+s_{f})x+1)^{2}}
\end{equation*}%
and the equation  $f^{\prime }(x)=0$ has a unique positive real solution
given by%
\begin{equation*}
x_{M}=\frac{1}{\sqrt{s_{h}}}.
\end{equation*}
As $s_{h}\leq 1$, we have that $\sqrt{s_{h}}\leq 1$ and so $x_{M}\geq 1$. As
a consequence, $f$ is strictly increasing in $[0,1]$. Then%
\begin{equation*}
\max_{x\in \lbrack 0,1]}f(x)=f(1)=1-\mu \leq 1.
\end{equation*}%
Then $f([0,1])\subseteq \lbrack 0,1]$. The equality is obtained when $\mu =0$.

\item[(F3)] The fixed points of $f$ are $x=0$ and%
\begin{equation*}
\overline{x}_{\pm }(\mu ,s_{f},s_{h})=\frac{s_{h}+s_{f}\pm \sqrt{%
(s_{h}-s_{f})^{2}-4s_{h}\mu (1-s_{f})}}{2s_{h}},
\end{equation*}%
which are reals if%
\begin{equation*}
\frac{(s_{h}-s_{f})^{2}}{4s_{h}(1-s_{f})}=\mu ^{\ast }\geq \mu .
\end{equation*}

\begin{enumerate}
\item If $\mu >\mu ^{\ast }$, then the map $f$ has two fixed points $x=0$
which is attracting and $x=1$ which is repelling.

\item If $\mu =\mu ^{\ast }$, then there exists one fixed point%
\begin{equation*}
\overline{x}(\mu ^{\ast },s_{f},s_{h})=\frac{s_{h}+s_{f}}{2s_{h}}=\frac{1}{2}%
+\frac{s_{f}}{2s_{h}}\geq \frac{1}{2}.
\end{equation*}%
In addition, $\overline{x}(\mu ^{\ast },s_{f},s_{h})\in (0,1)$ if and only
if $s_{f}<s_{h}$. Otherwise $\overline{x}(\mu ^{\ast },s_{f},s_{h})\geq 1$.

\item If $\mu <\mu ^{\ast }$, then the fixed points $\overline{x}_{-}(\mu
,s_{f},s_{h})$ and $\overline{x}_{+}(\mu ,s_{f},s_{h})$ belong to $(0,1)$ if
and only if $s_{f}<s_{h}$.
\end{enumerate}

\item[(F4)] When $\mu =0$, the fixed points are
\begin{equation*}
\overline{x}_{\pm }(0,s_{f},s_{h})=\frac{s_{h}+s_{f}\pm |s_{h}-s_{f}|}{2s_{h}%
}.
\end{equation*}%
If $s_{f}<s_{h}$, then $\overline{x}_{-}(0,s_{f},s_{h})=\frac{s_{f}}{s_{h}}<1
$ and $\overline{x}_{+}(0,s_{f},s_{h})=1$. Otherwise $\overline{x}%
_{-}(0,s_{f},s_{h})=1$ and $\overline{x}_{+}(0,s_{f},s_{h})>1$.
\end{enumerate}

The dynamics of $f$ is then simple. The fixed point $0$ satisfies that%
\begin{equation*}
1>f^{\prime }(0)=(\mu -1)(s_{f}-1)>0,
\end{equation*}%
and so it is always attracting. If there are no additional fixed points in $(0,1)$, then $0$ attracts all the orbits in $(0,1)$. If $f(1)\neq 1$, then
the orbit of $1$ also converges to $0$. If we have a unique non-zero fixed
point $\overline{x}(\mu ^{\ast },s_{f},s_{h})$ in $(0,1]$, then $\overline{x}(\mu ^{\ast },s_{f},s_{h})$ follows a fold bifurcation (see \cite{kutne}).

In \cite{zheng}, a new non-autonomous model is introduced by assuming that the parameters $\mu $, $s_{f}$ and $s_{h}$ can vary with time, that is, we have sequences $\mu _{n}$, $s_{f_{n}}$ and $s_{h_{n}}$ so that the model is given by a non-autonomous difference equation
\begin{equation}
x_{n+1}=f_{n}(x_{n})=\frac{(1-\mu _{n})(1-s_{f_{n}})x_{n}}{s_{h_{n}}x_{n}^{2}-(s_{h_{n}}+s_{f_{n}})x_{n}+1}.  \label{model1}
\end{equation}
When the sequences $\mu _{n}$, $s_{f_{n}}$ and $s_{h_{n}}$ are periodic with period $T$, the model is driven by a periodic sequence of maps $(f_{1},...f_{T},f_{1},...f_{T},...)$, denoted $[f_{1},...,f_{T}]$ from now
on. Several authors have studied periodic non-autonomous models in biology \cite{ziyad3,ziyad2,Canovas1,Canovas3,elaysac1,elaysac2,elaysac3,Luis,mendoza,silva}. Regarding the model in \cite{zheng}, the dynamics of the non-autonomous model are more complicated, and the following conjectures were stated.

\begin{conjecture}
\label{conjetura}Assume that%
\begin{equation*}
0\leq s_{f_{n}}<s_{h_{n}}\leq 1\ \ \ \text{\ and\ \ \ }\mu _{n}\leq \mu _{n}^{\ast }=\frac{(s_{h_{n}}-s_{f_{n}})^{2}}{4s_{h_{n}}(1-s_{f_{n}})}
\end{equation*}%
and $\mu _{n}$, $s_{f_{n}}$ and $s_{h_{n}}$ are periodic with period $T$.
Then model (\ref{model1}) has at most two $T$-periodic solutions. To be
precise, if $\mu_{n}\neq 0$, then model (\ref{model1}) has exactly two $T$-periodic
solutions, among which one is stable and the other is unstable. If $\mu_{n}=0$, then model (\ref{model1}) has a unique $T$-periodic solution, which
is unstable.
\end{conjecture}

In this paper, we disprove this conjecture. We also study how to solve the
problems in this conjecture by analyzing the conditions that a
periodic model (\ref{model1}) has to satisfy to show a scenario slightly
similar to that proposed in the conjecture. Basic notation, known results
and technical lemmas are given in the next section, while the last section
is devoted to proving our main results and additional discussion.

\section{Basic notation, known results and technical lemmas}

Given a non-autonomous difference
equation $x_{n+1}=\varphi _{n}(x_{n})$, with $\varphi _{n}:I=[a,b]\subset
\mathbb{R}\rightarrow \lbrack a,b]$, we define the orbit of $x\in \lbrack
a,b]$ as the sequence $(x_{n})$ given by $x_{1}=x$ and $x_{n+1}=\varphi
_{n}(x_{n})$ for $n\geq 1$. The set of limit points of the orbit is the $\omega $-limit set, denoted $\omega (x,(\varphi _{n})).$ An orbit is periodic if
there is a minimal positive integer $T$ such that $x_{n+T}=x_{n}$ for all $%
n\geq 1$. In this case, $T$ is the period of $(x_{n})$. When the period $T=1$, the sequence is constant, and we say that $x$ is a fixed point. When
the sequence $(\varphi _{n})$ is constant, that is $\varphi_n=\varphi$ for $n\in \mathbb{N}$, we obtain analogous definitions for
an autonomous difference equation $x_{n+1}=\varphi (x_{n}).$

Following~\cite{Milnor}, we define \emph{(minimal) metric attractor} as
follows.

\begin{definition}
Let $\varphi :I\rightarrow I$ be a continuous map. A forward invariant
closed subset $A\subset I$ is said to be a \emph{metric attractor} if the
following two conditions hold:

\begin{enumerate}
\item[(a)] The basin of attraction of $A$, namely $B(A):=\{x:\,\omega
(x,\varphi )\subset A\}$, has positive Lebesgue measure;

\item[(b)] For every forward invariant compact set $A^{\prime }\subset A$
strictly contained in $A$ we have that $B(A^{\prime })$ has Lebesgue
measure equal to zero.
\end{enumerate}
\end{definition}

It is well-known that a fixed point $x$ of a $C^{1}$ interval map $\varphi $
such that $|\varphi ^{\prime }(x)|<1$ is an attractor of $\varphi $ (see
e.g. \cite{Elaydi0}), and we say that $x$ is attracting. If $|\varphi
^{\prime }(x)|>1$ we say that the fixed point is repelling. For
non-hyperbolic fixed points, such that $|\varphi ^{\prime }(x)|=1$
generalized conditions to guarantee that $x$ is an attractor can be also
found in \cite{Elaydi0}. If the point $x$ is periodic with period $T$, this
conditions are replaced by $|(\varphi ^{n})^{\prime }(x)|<1$ and $|(\varphi
^{n})^{\prime }(x)|>1$ in the hyperbolic case, where $\varphi ^{n}$ is the
$n$-th composition of $\varphi $.

The sequence $(\varphi _{n})$ is periodic if there exists a minimal positive
integer $T$ such $\varphi _{n+T}=\varphi _{n}$ for all $n\geq 1$. We denote
periodic sequences of maps as $[\varphi _{1},...,\varphi _{T}]$. The
following result can be found in \cite{ziyad} or \cite{canovas}.

\begin{proposition}
\label{periodic}Let $[\varphi _{1},...,\varphi _{T}]$ be a periodic sequence
of maps and let $x$ be a fixed point of $\varphi _{T}\circ ...\circ \varphi
_{1}$. Then, $x$ is a fixed point for all the maps $\varphi _{1},...,\varphi
_{T}$, if and only if $x$ is a fixed point for $[\varphi _{1},...,\varphi
_{T}]$.
\end{proposition}

The notion of Schwarzian derivative was used by D. Singer~\cite{Singer} to
get information about the number of attracting periodic orbits.

\begin{definition}
Let $\varphi :I\rightarrow I$ be a $\mathcal{C}^{3}$ map defined on an
interval $I$, then
\begin{equation*}
\mathcal{S}(\varphi )(x)=\frac{\varphi ^{\prime \prime \prime }(x)}{\varphi
^{\prime }(x)}-\frac{3}{2}\left( \frac{\varphi ^{\prime \prime }(x)}{\varphi
^{\prime }(x)}\right) ^{2},
\end{equation*}%
if $\varphi ^{\prime }(x)\neq 0$. Moreover, $\varphi $ is said to have
negative Schwarzian derivative, denoted by $\mathcal{S}(\varphi )<0$, if $%
\mathcal{S}(\varphi )(x)<0$ for each $x\in I$ such that $\varphi ^{\prime
}(x)\neq 0$.
\end{definition}

The composition of maps with negative Schwarzian derivative also has its
Schwarzian derivative negative since it is held that
\begin{equation}
\mathcal{S}(\varphi _{2}\circ \varphi _{1})=(\varphi _{1}^{\prime })^{2}(x)\mathcal{S}(\varphi _{2})(\varphi _{1}(x))+\mathcal{S}(\varphi _{1})(x).
\label{comder}
\end{equation}%
for $C^{3}$ interval maps $\varphi _{1}$ and $\varphi _{2}$.

A map $\varphi :I\rightarrow I$, $I=[a,b]$, $a,b\in \mathbb{R}$, is said to be unimodal if there exist $c\in (a,b)$ such that $\varphi$ is increasing (resp. decreasing) in $[a,c)$ and decreasing (resp. increasing) in $(c,b]$. Note that if $\varphi$ is $C^1$, we have that $\varphi ^{\prime}(c)=0$.

\begin{theorem}[Singer's Theorem]
Let $I$ be a closed interval and $\varphi :I\rightarrow I$ be a $\mathcal{C}^{3}$ unimodal map such that $\mathcal{S}(\varphi )<0$. Then each
attracting periodic orbit attracts at least the one critical point or a
boundary point.
\end{theorem}

As a consequence of Singer's Theorem, when $\varphi $ is a $\mathcal{C}^{3}$
unimodal map with $\mathcal{S}(\varphi )<0$, the map $\varphi $ can have at most two attracting fixed points, and one of them belongs to the boundary of $I$.

\begin{lemma}
\label{unimod1}Let $\varphi :[a,b]\rightarrow \lbrack a,b]$ be unimodal with
maximum point $c$ and negative Schwarzian derivative. Assume that $a$ is an
attracting fixed point of $\varphi $ and $\varphi (c)\leq c$. Then $\varphi $
has at most three fixed points smaller than $c$.
\end{lemma}

\textbf{Proof.} By Singer's Theorem, the map $\varphi $ has at most two
metric attractors. Clearly, the fixed point $a$ is one attractor. As $\varphi
(c)\leq c$, all the fixed points of $\varphi $ are located in $[a,c]$
because $\varphi (x)<x$ for all $x\in (c,b]$. Note that $\varphi $ increases on $[a,c]$. Assume that $a<x_{1}<x_{2}<x_{3}$ are fixed
points of $\varphi $ with $x_{3}\leq c$. We distinguish several cases.

\begin{enumerate}
\item[(a)] Assume that $\varphi (x)\leq x$ for all $x\in \lbrack a,c]$. Then
any orbit starting in $(a,x_{1})$ converges to $a$, any orbit starting in $%
(x_{1},x_{2})$ converges to $x_{1}$ and any orbit starting in $(x_{2},x_{3})$
converges to $x_{2}$. Thus, $\varphi $ would have three metric attractors,
which is impossible.

\item[(b)] So, there is $x^{\ast }>x_{1}$ such that $\varphi (x^{\ast
})>x^{\ast }$. If $x^{\ast }<x_{2}$ and assume that $\varphi (x)\geq x$ for
all $x\in (x_{1},c)$, then any orbit starting in $(x_{1},x_{2})$ converges
to $x_{2}$ and any orbit starting in $(x_{2},x_{3})$ converges to $x_{3}$,
and again $\varphi $ would have three metric attractors. So there must exist
$x_{\ast }>x_{2}$ such that $\varphi (x_{\ast })<x_{\ast }$. Here we have
two cases. If $\varphi (x)\leq x$ for all $x\in (x_{2},c)$, then $x_{3}\neq c
$ because $\varphi ^{\prime }(x_{3})=1$ and hence any orbit starting in $%
(x_{1},x_{3})$ would converge to $x_{2}$ and any orbit in $(x_{3},c)$ would
converge to $x_{3}$ and again we would have three attractors. If there
exists $x_{\ast }^{\ast }>x_{3}$ such that $\varphi (x_{\ast }^{\ast
})>x_{\ast }^{\ast }$, then as $\varphi (c)\leq c$ we would have another
fixed point $x_{4}\in (x_{\ast }^{\ast },c]$, and hence $x_{2}$ and $x_{4}$
would be metric attractors, which is impossible.

\item[(c)] If $x^{\ast }\in (x_{2},x_{3})$, then $x_{1}$ attracts the orbits
starting in $(x_{1},x_{2})$ and $x_{3}$ attracts the orbits starting in $%
(x_{2},x_{3})$, again we would have three metric attractors, which is
impossible.

\item[(d)] Finally, if $x^{\ast }>x_{3}$, then $x_{1}$ attracts the
orbits starting in $(x_{1},x_{2})$ and $x_{2}$ attracts the orbits starting
in $(x_{2},x_{3})$, again we would have three metric attractors, which is
impossible.
\end{enumerate}

Since the above cases cover all the possible scenarios, the proof concludes.$%
\square $

Next, we introduce some basic lemmas on the models we are dealing with.

\begin{lemma}
\label{schwarz1}Let $f$ be the map given in (\ref{map}). Then $\mathcal{S}(f)(x)<0$ whenever $f^{\prime}(x)\neq 0$.
\end{lemma}

\textbf{Proof.} A direct computation gives us that%
\begin{equation*}
\mathcal{S}(f)(x)=\frac{f^{\prime \prime \prime}(x)}{f^{\prime }(x)}-\frac{3}{2}\left( \frac{%
f^{\prime \prime }(x)}{f^{\prime }(x)}\right) ^{2}=-\frac{6s_{h}}{%
(s_{h}x^{2}-1)^{2}}<0,
\end{equation*}%
for all $x\neq 1/\sqrt{s_{h}}$.$\square $

\begin{lemma}
\label{schwarz2}Let $[f_{1},...,f_{T}]$ be the periodic model given in (\ref%
{model1}). Then $\mathcal{S}(f_{T}\circ ...\circ f_{1})<0$.
\end{lemma}

\textbf{Proof.} It is a consequence of Lemma \ref{schwarz1} and \ (\ref{comder}).$\square $

\begin{lemma}
Let $f$ be the map (\ref{map}) and assume that $0\leq s_{f}<s_{h}\leq 1$. Let $x_{M}=%
\frac{1}{\sqrt{s_{h}}}$ be its positive maximum. Then%
\begin{equation*}
f\left( \frac{1}{\sqrt{s_{h}}}\right) \leq \frac{1}{\sqrt{s_{h}}}.
\end{equation*}
\end{lemma}

\textbf{Proof.} Note that,
\begin{equation*}
f\left( \frac{1}{\sqrt{s_{h}}}\right) =\frac{(1-\mu )(1-s_{f})}{2\sqrt{s_{h}}%
-(s_{h}+s_{f})}.
\end{equation*}%
To check the inequality, note that%
\begin{equation*}
\frac{(1-\mu )(1-s_{f})}{2\sqrt{s_{h}}-(s_{h}+s_{f})}\leq \frac{(1-s_{f})}{2%
\sqrt{s_{h}}-(s_{h}+s_{f})} .
\end{equation*}%
We have to prove that%
\begin{equation*}
\frac{(1-s_{f})}{2\sqrt{s_{h}}-(s_{h}+s_{f})}\leq \frac{1}{\sqrt{s_{h}}} ,
\end{equation*}%
which is true if and only if%
\begin{equation*}
s_{f}(1-\sqrt{s_{h}})\leq \sqrt{s_{h}}-s_{h},
\end{equation*}%
that is%
\begin{equation*}
s_{f}\leq \frac{\sqrt{s_{h}}-s_{h}}{1-\sqrt{s_{h}}}=\sqrt{s_{h}},
\end{equation*}%
which is true because $s_{f}<s_{h}\leq \sqrt{s_{h}}$.$\square $

We prove our main result in the next section.

\section{Main results}

We consider periodic sequences $\mu _{n}$, $s_{f_{n}}$ and $s_{h_{n}}$ such
that there exists $T$ with the property that%
\begin{equation}
\mu _{n}=\mu _{n+T},\ \ s_{f_{n}}=s_{f_{n+T}},\ \ s_{h_{n}}=s_{h_{n+T}}.
\label{cond1}
\end{equation}%
We denote by $f_{n}$ the model with parameters $\mu _{n}$, $s_{f_{n}}$ and $s_{h_{n}}$, which clearly is a periodic sequence of maps of period $T$. Let us also assume the hypothesis of Conjecture \ref{conjetura}. The first example disproves the $\mu _{n}\neq 0$ case.

\begin{example}
Fix $T=2$, $s_{f_{1}}=s_{f_{2}}=\frac{1}{20}$, $s_{h_{1}}=\frac{9}{10}$, $s_{h_{2}}=\frac{3}{10}$, $%
\mu _{1}=\frac{(s_{h_{1}}-s_{f_{1}})^{2}}{4s_{h_{1}}(1-s_{f_{1}})}$ and $\mu
_{2}=\frac{(s_{h_{2}}-s_{f_{2}})^{2}}{4s_{h_{2}}(1-s_{f_{2}})}$. It is easy
to see that the unique real solution of the equation $(f_{2}\circ
f_{1})(x)=x$ is $x=0$. There are four complex solutions given by the roots of the polynomial
$$
-4523020 + 21055109 x - 34761128 x^2 + 26901936 x^3 - 11197440 x^4.
$$
They are two pairs of
complex conjugate solutions, given approximately by $0.539661 \pm 0.0228932 i$ and $0.661593 \pm 0.973024 i$. Thus, the fixed point
$0$ is the global attractor of this system.

Additionally, fix $T=2$, $s_{f_{1}}=s_{f_{2}}=\frac{1}{20}$, $s_{h_{1}}=\frac{9}{10}$, $%
s_{h_{2}}=\frac{3}{10}$, $\mu _{1}=\frac{(s_{h_{1}}-s_{f_{1}})^{2}}{%
4s_{h_{1}}(1-s_{f_{1}})}-10^{-9}$ and $\mu _{2}=\frac{%
(s_{h_{2}}-s_{f_{2}})^{2}}{4s_{h_{2}}(1-s_{f_{2}})}$. Note that $\mu_1<\mu_1^*$ and again, $0$ is the unique
solution of the equation $(f_{2}\circ f_{1})(x)=x$.
\end{example}

We prove the next result.

\begin{theorem}
\label{main1}Assume $T=2$ and the hypothesis of Conjecture \ref{conjetura}
such that $\mu _{1}=\mu _{2}=0$. Then $f_{2}\circ f_{1}$ has one fixed point
in $(0,1)$.
\end{theorem}

\textbf{Proof.} First, note that $1$ is a fixed point for both maps $f_{1}$
and $f_{2}$. As%
\begin{equation*}
0<(f_{2}\circ f_{1})^{\prime }(1)=f_{2}^{\prime }(1)f_{1}^{\prime }(1)=\frac{%
1-s_{h_{2}}}{1-s_{f2}}\frac{1-s_{h_{1}}}{1-s_{f_{1}}}<1,
\end{equation*}%
the fixed point $1$ must be attracting. This would be impossible if $%
(f_{2}\circ f_{1})(x)<x$ for all $x\in (0,1)$. As a consequence, there must
be $x^{\ast }\in (0,1)$ such that $(f_{2}\circ f_{1})(x^{\ast })=x^{\ast }$.
We have to check that $x^{\ast }$ is the unique fixed point in $(0,1)$.

First, assume that $s_{h_{2}}\leq s_{h_{1}}$, then%
\begin{equation*}
\frac{1}{\sqrt{s_{h_{1}}}}\leq \frac{1}{\sqrt{s_{h_{2}}}},
\end{equation*}%
so, the equation%
\begin{equation*}
(f_{2}\circ f_{1})^{\prime }(x)=f_{2}^{\prime }(f_{1}(x))f_{1}^{\prime }(x)=0
\end{equation*}%
has a unique solution given by
\begin{equation*}
f_{1}^{\prime }(x)=0
\end{equation*}%
because $f_{1}(x)<\frac{1}{\sqrt{s_{h_{2}}}}$ for all $x$. Then, $f_{2}\circ
f_{1}$, defined in $[0,z]$, where $z>\frac{1}{\sqrt{s_{h_{1}}}}$ is such that $(f_{2}\circ f_{1})([0,z])\subset [0,z]$,  is unimodal and has negative Schwarzian
derivative. By Lemma \ref{unimod1}, $f_{2}\circ f_{1}$ can have at most
three fixed points smaller than $\frac{1}{\sqrt{s_{h_{1}}}}>0$. So, these
fixed points must be $0,\ x^{\ast }$ and $1$. As the number fixed points of $%
f_{2}\circ f_{1}$ and $f_{1}\circ f_{2}$ are the same, the case $%
s_{h_{2}}\geq s_{h_{1}}$ follows taking $f_{1}\circ f_{2}$ and the proof
concludes.$\square $

\begin{figure}[htbp]
\begin{center}

\includegraphics[width=0.35\textwidth]{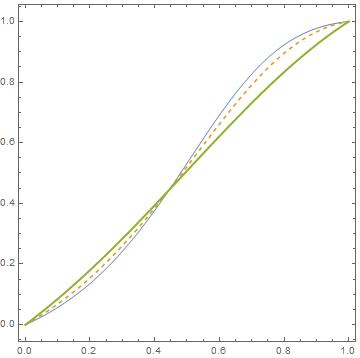}

\caption{For $T=2$, $s_{f_{1}}=0.2$, $s_{h_{1}}=0.45$, $s_{f_{2}}=0.4$, $s_{h_{2}}=0.9$ and $\mu_1=\mu_2=0$, graph of the maps $f_1$ (thick), $f_2$ (dashed) and $f_2\circ f_1$ (thin). The common fixed point is $4/9$.
 }
\label{fig1}
\end{center}
\end{figure}

\begin{example}
\label{ex2}Fix $T=2$, $\mu _{1}=\mu _{2}=0$ and the parameter fulfilling the
conditions of Conjecture \ref{conjetura}. The map $f_{2}\circ f_{1}$ has a
unique fixed point $x^{\ast }\in (0,1)$. If $f_{2}(x)\neq f_{1}(x)$, this
point generates a periodic orbit of period two of $[f_{1},f_{2}]$. However,
is it possible that $f_{2}(x^{\ast })=f_{1}(x^{\ast })=x^{\ast }$.  We
consider the solutions%
\begin{equation*}
\overline{x}_{-}(0,s_{f_{i}},s_{h_{i}})=\frac{s_{f_{i}}}{s_{h_{i}}},\ i=1,2.
\end{equation*}%
Assume that%
\begin{equation}
\frac{s_{f_{1}}}{s_{h_{1}}}=\frac{s_{f_{2}}}{s_{h_{2}}}\text{,}
\label{cond2}
\end{equation}%
Then the fixed point of $f_{2}\circ f_{1}$ is $x^{\ast }=\overline{x}_{-}(0,s_{f_{i}},s_{h_{i}})$ and generates a repelling fixed point of $[f_{1},f_{2}]$. So, the second part of Conjecture \ref{conjetura} is also
false, although one repelling periodic orbit exists with a period of either 1 or 2. For instance, for $s_{f_{1}}=0.2$, $s_{h_{1}}=0.45$, $s_{f_{2}}=0.4$ and $s_{h_{2}}=0.9$, the common fixed point is repelling, as Figure \ref{fig1} shows.
\end{example}

\begin{remark}
Example \ref{ex2} shows that the statement of Theorem 2.3 from \cite{zheng} should be more precise. The reader could think that the periodic solution is of period two, but this repelling periodic solution can be, in fact, a repelling fixed point.
\end{remark}

Following the above idea, we can consider the case where $\mu _{1}=0$ and $%
\mu _{2}=\mu _{2}^{\ast }$. Then, the unique fixed point of $f_{1}$ is%
\begin{equation*}
x_{1}=\overline{x}_{-}(0,s_{f_{1}},s_{h_{1}})=\frac{s_{f_{1}}}{s_{h_{1}}}
\end{equation*}%
and the non-zero fixed points of $f_{2}$ is%
\begin{equation*}
x_{2}=\overline{x}(\mu _{2}^{\ast },s_{f_{2}},s_{h_{2}})=\frac{%
s_{h_{2}}+s_{f_{2}}}{2s_{h_{2}}}.
\end{equation*}%
If we impose the condition $x_{1}=x_{2}$, we have that%
\begin{equation*}
\frac{s_{f_{1}}}{s_{h_{1}}}=\frac{s_{h_{2}}+s_{f_{2}}}{2s_{h_{2}}},
\end{equation*}%
and this condition implies that a fixed point of $f_{2}\circ f_{1}$ is also
a common fixed of $f_{1}$ and $f_{2}$, while there is another fixed point of
$f_{2}\circ f_{1}$. The first fixed point generates a fixed point of $%
[f_{1},f_{2}]$ while the second one generates a periodic orbit of period
two. For instance, for $s_{f_{1}}=0.5$, $s_{h_{1}}=0.8$, $s_{f_{2}}=0.2$ and
$s_{h_{2}}=0.8$, we have that the common fixed point is repelling and the
periodic point is attracting.

Similar conditions can be posed to the other cases when $\mu _{1}$ and $\mu
_{2}$ are positive. So, we can find parameter values for which the fixed
points of $f_{2}\circ f_{1}$ produces either a fixed point or a periodic
orbit of period two. In addition, we can prove the following result.

\begin{theorem}\label{mainth}
Let $T=2$ and let the hypothesis of Conjecture \ref{conjetura} be satisfied
with $\mu _{i}\neq 0$, $i=1,2$. Then, the number of non-zero fixed points of
$f_{2}\circ f_{1}$ is at most two.
\end{theorem}

{\bf Proof.} The proof is analogous to that of Theorem \ref{main1}, with one extra case. Here we can have none, one or two fixed points. The bound of the number of fixed points follows the same argument than that of Theorem \ref{main1}.
$\square$

\begin{remark}
If $\mu _{1}\neq 0$ and $\mu _{2}=0$, if%
\begin{equation*}
\frac{s_{f_{2}}}{s_{h_{2}}}>1-\mu _{1},
\end{equation*}%
then $f_{2}\circ f_{1}$ does not have fixed points different from zero. Note that $f_1([0,1])=[0,1-\mu_1]$ and $f_2(x)<x$ for all $x\in (0,\frac{s_{f_{2}}}{s_{h_{2}}})$, and so $(f_2\circ f_1)(x)<x$ for all $x\in (0,1]$. An example of this case can be seen in Figure \ref{fig2}(a).

\begin{figure}[htbp]
\begin{center}
\begin{tabular}{cc}
(a) \includegraphics[width=0.35\textwidth]{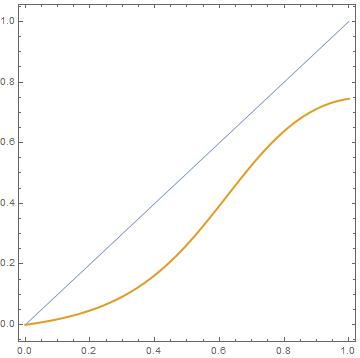} &(b) \includegraphics[width=0.35\textwidth]{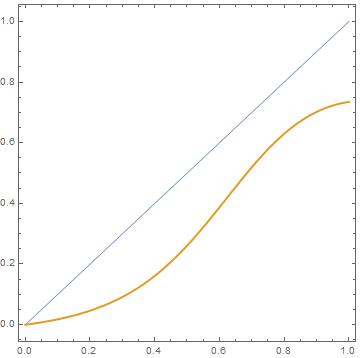}
\end{tabular}
\caption{(a) For $T=2$, $s_{f_{1}}=0.1$, $s_{h_{1}}=0.9$, $s_{f_{2}}=0.3$, $s_{h_{2}}=0.9$, $\mu_1=\mu_1^*$ and $ \mu_2=0$, graph of the maps $f_2\circ f_1$ (thick) and the identity (thin).
(b) For $T=2$, $s_{f_{1}}=0.1$, $s_{h_{1}}=0.9$, $s_{f_{2}}=0.3$, $s_{h_{2}}=0.9$, $\mu_1=\mu_1^*$ and $ \mu_2=\mu_2^*$, graph of the maps $f_2\circ f_1$ (thick) and the identity (thin). There are no fixed points except for $0$.
 }
\label{fig2}
\end{center}
\end{figure}

If $%
\mu _{2}\neq 0$, the condition%
\begin{equation*}
\overline{x}_{-}(\mu _{2},s_{f_{2}},s_{h_{2}})>1-\mu _{1}
\end{equation*}%
also implies that $f_{2}\circ f_{1}$ does not have fixed points different
from zero. In particular, if $\mu _{2}=\mu _{2}^{\ast }$, then the condition
reads as%
\begin{equation*}
\frac{s_{h_{2}}+s_{f_{2}}}{2s_{h_{2}}}>1-\mu _{1}.
\end{equation*}%
An example of this case can be seen in Figure \ref{fig2}(a).
If%
\begin{equation*}
\overline{x}_{-}(\mu _{2},s_{f_{2}},s_{h_{2}})\leq 1-\mu _{1}<\overline{x}%
_{+}(\mu _{2},s_{f_{2}},s_{h_{2}})
\end{equation*}%
we could have just two fixed points. In that case, the non-zero fixed point
must have slope one and $(f_{2}\circ f_{1})(x)\leq x$ for all $x\in \lbrack
0,1]$. For instance, for $s_{f_1} = 0.1$, $s_{h_1} = 0.9$, $s_{f_2} = 0.8$, $s_{h_2} = 0.9$, $\mu_1 =0.0975309$ and $\mu_2=0.00863972$ we have a unique non-zero fixed point $x\approx 0.7949203$. Figure \ref{fig4} shows the graph of the map $f_2\circ f_1$ in this case.

\begin{figure}[htbp]
\begin{center}

\includegraphics[width=0.35\textwidth]{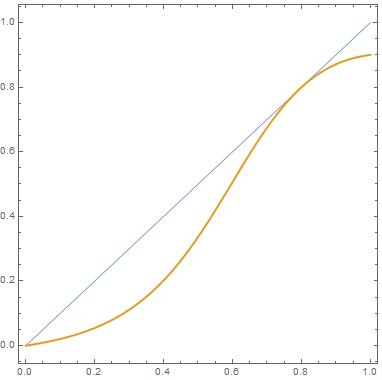}

\caption{For $T=2$, $s_{f_1} = 0.1$, $s_{h_1} = 0.9$, $s_{f_2} = 0.8$, $s_{h_2} = 0.9$, $\mu_1 =0.0975309$ and $\mu_2=0.00863972$, graph of the map $f_2\circ f_1$ (thick) and the identity (thin). The fixed point is approximately $0.7949203$.
 }
\label{fig4}
\end{center}
\end{figure}

\end{remark}

Finally, for $T\geq 3$, we can prove the following result, whose proof is
analogous to that of Theorem \ref{main1}.


\begin{theorem}
Let $T\geq 3$ and let the hypothesis of Conjecture \ref{conjetura} be
satisfied. Then the map
$f_{T}\circ ...\circ f_{1}$ has at most two non-zero fixed points in $(0,1]$.
\end{theorem}

{\bf Proof.}
By \cite{zheng}, all the fixed points of $f_T\circ ...\circ f_1$ must be contained in $[0,1]$. Let $x_m>1$ be the first extremum of $f_T\circ ...\circ f_1$. Clearly $x_m$ is a maximum and $(f_T\circ ...\circ f_1)(x_m)<x_m$, because otherwise we would have a fixed point greater than one. Let $z>x_m$ be such that the remaining extrema of the map are greater than $z$. Then, $(f_T\circ ...\circ f_1)([0,z])\subset [0,z]$ and $(f_T\circ ...\circ f_1)|_{[0,z]}$ is unimodal. Following the proof of Theorem \ref{main1}, we conclude that the map has at most two non-zero fixed points.
$\square$

\begin{remark}
We want to thank the anonymous referee who showed the existence of reference \cite{sacker}. This paper can help answer the following question. We have shown that the periodic model can have the fixed point $0$ as a global attractor. However, the conditions on the parameters need to be clarified to guarantee this fact. This problem deserves further investigation.
\end{remark}

\section{Discussion}

We have shown that Conjecture \ref{conjetura} is not true and proved alternative results. It must be emphasized that given a periodic difference equation given by $[f_1,...,f_T]$, its dynamics is related to the composition $f_T\circ ...\circ f_1$ and it is quite difficult to characterize it in terms of the individual maps $f_1,...,f_T$. It is also important to remark that the technical difficulties in analyzing the dynamics of $f_T\circ ...\circ f_1$ increase when $T$ increases. So, analyzing periodic difference equations, depending on several parameters, generally becomes quite challenging.

\section*{Acknowledgements}

We thank the anonymous referees of this paper for their suggestions and remarks that helped us polish the final paper version.

J.S.C\'anovas' address: {Departamento de Matem\'atica Aplicada y Estad\'{i}stica, Universidad Poli\-t\'ec\-nica de Cartagena, C/ Doctor Fleming sn, 30202, Cartagena, Spain. Email:
Jose.Canovas@upct.es}
\end{document}